\documentclass[12pt]{article}
\usepackage{amsfonts}
\usepackage{amssymb}
\usepackage{amsmath}
\usepackage{eufrak}
\usepackage[dvips]{graphicx}
\usepackage{latexsym}
\newcommand{\ind}{\makebox[1em]{\raisebox{-.5ex}[0ex][0ex]{\makebox[0em]%
{$\smile$}}\raisebox{.4ex}[0ex][0ex]{\makebox[-.02em]{$|$}}}}

\newcommand{\dep}{\makebox[1em]{\raisebox{.3ex}[0ex][0ex]%
{$\not$}\makebox[.7em]{\ind}}}

\newcommand{\nmdep}{\makebox[1em]{\raisebox{1.5ex}[0ex][0ex]{\makebox[0em]%
{$\scriptscriptstyle n\! m$}}\makebox[-1em]{$\dep$}}}
\newcommand{\nmind}{\makebox[1em]{\raisebox{1.5ex}[0ex][0ex]{\makebox[0em]%
{$\scriptscriptstyle n\! m$}}\makebox[-1em]{$\ind$}}}

\newcommand{\bi}{\begin{itemize}}
\newcommand{\ei}{\end{itemize}}

\newcommand\NM{\mathcal{N\!M}}
\newtheorem{theorem}{Theorem}[section]
\newtheorem{lemma}[theorem]{Lemma}
\newtheorem{fact}[theorem]{Fact}

\newtheorem{corollary}[theorem]{Corollary}
\newtheorem{proposition}[theorem]{Proposition}
\newtheorem{definition}[theorem]{Definition}
\newtheorem{remark}[theorem]{Remark}
\newtheorem{conjecture}[theorem]{Conjecture}

\def\defn{\begin{definition}\upshape}
\def\edefn{\end{definition}}

\title{Small, $nm$-stable compact $G$-groups}
\author{Krzysztof Krupi\'nski\footnote{Research supported by the Polish Government grant N N201 545938}  \hspace{0.1pt} and Frank Wagner\footnote{Research partially supported by European Community Network MRTN-CT-2004-512234 Modnet and the French Research Agency grant ANR-09-BLAN-0047 Modig}}
\date{}

\begin{document}
\maketitle
\begin{abstract}
We prove that if $(H,G)$ is a small, $nm$-stable compact $G$-group, then $H$ is nilpotent-by-finite, and if additionally $\NM(H) \leq \omega$, then $H$ is abelian-by-finite.
Both results are significant steps towards the proof of the conjecture that each small, $nm$-stable compact $G$-group is abelian-by-finite.

We provide counter-examples to the $\NM$-gap conjecture, that is we give examples of small, $nm$-stable compact $G$-groups of infinite ordinal $\NM$-rank.
\end{abstract}

\section[\mbox{}]{Introduction}

In \cite{kr}, the notion of a Polish structure and a ternary relation of independence (called $nm$-independence) were introduced in order to apply model-theoretic ideas to purely topological and descriptive set theoretic objects, such as Polish $G$-spaces. Although we shall review the basic definitions and results, a certain familiarity with \cite{kr} will be helpful.

\defn
A {\em Polish structure} is a pair $(X,G)$, where $G$ is a Polish group acting faithfully on a set $X$ so that the stabilizers of all singletons are closed subgroups of $G$. We say that $(X,G)$ is {\em small} if for every $n \in \omega$ there are only countably many orbits on $X^n$ under the action of $G$. 
\edefn

The class of Polish structures is much wider than the class of profinite structures introduced by Newelski \cite{ne01,ne02} and contains many natural examples \cite{kr,ca}. Typical examples are of the form $(X,Homeo(X))$, where $X$ is a compact metric space and $Homeo(X)$ is its group of homeomorphisms equipped with the compact-open topology, e.g. $(S^1, Homeo(S^1))$ or $(I^\omega, Homeo(I^\omega))$.

For a Polish structure $(X,G)$, $A\subseteq X$ and $a \in X^n$, $G_A$ denotes the pointwise stabilizer of $A$ in $G$ and $o(a/A):=\{ ga: g \in G_A\}$ is the orbit of $a$ over $A$; $o(a)$ stands for $o(a/\emptyset)$.

\defn
Let $(X,G)$ be a Polish structure, $a$ be a finite tuple and $A,B$ finite subsets of $X$. Let $\pi_A:G_A \to o(a/A)$ be defined by $\pi_A(g)=ga$. We say that $a$ is {\em $nm$-independent} from $B$ over $A$ (written $a\nmind_{A}B$) if $\pi_{A}^{-1}[o(a/AB)]$ is non-meager in $\pi_{A}^{-1}[o(a/A)]$. Otherwise, we say that $a$ is {\em $nm$-dependent} on $B$ over $A$ (written $a\nmdep_{A}B$).
\edefn

By \cite{kr}, $\nmind$ shares several nice properties with  forking independence in stable theories, e.g. invariance under $G$, symmetry and transitivity. Assuming  smallness, $\nmind$ also satisfies the existence of independent extensions. In a small Polish structure $(X,G)$, this leads to a counterpart of SU-rank, called $\NM$-rank, having all the expected  properties, including the Lascar Inequalities. 

\defn
The {\em rank $\cal{NM}$} is the unique function from the collection of orbits over finite sets to the ordinals together with $\infty$, satisfying
$$\begin{array}{lll}
\NM(a/A) \geq \alpha +1 & \mbox{iff}  &\mbox{there is a finite set}\;\; B \supseteq A \;\; \mbox{such that} \\
& &a \nmdep_A B\;\; \mbox{and} \;\; \NM(a/B) \geq \alpha . \end{array}$$
The {\em $\NM$-rank} of $X$ is defined as the supremum of $\NM(x/\emptyset)$, $x \in X$.
\edefn

This, in turn, leads to a counterpart of superstability, called {\em $nm$-stability}. A small Polish structure $(X,G)$ is said to be {\em $nm$-stable} if $\NM(X) <\infty$.

In this paper, we are going to study the structure of compact groups 
in the context of small, $nm$-stable Polish structures. Let us recall the relevant definitions.

\defn
A {\em compact $G$-group} 
is a Polish structure $(H,G)$, where $G$ acts continuously and by automorphisms on a compact group $H$. 
\edefn

\defn
A {\em profinite group regarded as profinite structure} is a pair $(H,G)$ such that $H$ and $G$ are profinite groups and $G$ acts continuously on $H$ as a group of automorphisms.
\edefn

In this paper, compact spaces and topological groups are Hausdorff by definition, and a profinite group is the inverse limit of a countable inverse system of finite groups.
When the context is clear, we skip the phrase `regarded as profinite structure'. 

It was noticed in \cite{kr} that if $(H,G)$ is a small compact $G$-group, then the group $H$ is locally finite and hence profinite. 
However, $G$ is only Polish (not necessarily compact), which makes the class of small compact $G$-groups much wider than the class of small profinite groups (regarded as profinite structures).

Recall that originally profinite groups regarded as profinite structures were defined as pairs $(H,G)$, where $H$ is a profinite group and $G$ is a closed subgroup of the group of all automorphisms of $H$ preserving a distinguished inverse system defining $H$. In particular, there is basis of open neighborhoods of $e$ in $H$ consisting of clopen, normal, invariant subgroups (which is in general not the case for small compact $G$-groups).

\cite[Theorem 1.9]{kr05a} says that each abelian profinite group of finite exponent considered together with the group of all automorphisms respecting a distinguished inverse system indexed by $\omega$ is small. On the other hand, the main conjecture on small profinite groups \cite{ne01} is the following.

\begin{conjecture}\label{groups}
Each small profinite group is abelian-by-finite.
\end{conjecture}

In \cite{wa03}, the conjecture was proved assuming additionally $nm$-stability.
We also have the following three intermediate conjectures, which are wide open.

\begin{conjecture}
Let $(H,G)$ be a small profinite group. Then:
\begin{enumerate}
\item[(A)] $H$ is solvable-by-finite.
\item[(B)] If $H$ is solvable-by-finite, it is nilpotent-by-finite.
\item[(C)] If $H$ is nilpotent-by-finite, it is abelian-by-finite.
\end{enumerate}
\end{conjecture}

In this paper, we will consider generalizations of these conjectures to the wider class of small compact $G$-groups. 
It turns out that without any extra assumptions, these generalizations are false \cite{kr}. So, the question is whether they are true under the additional assumption of $nm$-stability. We will denote the conjectures obtained in this way by $(A')$, $(B')$ and $(C')$.

\cite[Theorem 5.19]{kr} says that $(A')$ is true, and \cite[Theorem 5.24]{kr} that $(B')$ is true assuming additionally that $\NM(H) <\omega$.

In this paper, we prove $(B')$ in its full generality and $(C')$ assuming additionally that $\NM(H) \leq \omega$. 
As a consequence, we obtain the following two theorems.

\begin{theorem}\label{Theorem X}
If $(H,G)$ is a small, $nm$-stable compact $G$-group, then $H$ is nilpotent-by-finite.
\end{theorem}

\begin{theorem}\label{Theorem XX}
If $(H,G)$ is a small compact $G$-group, and either $\NM(H) < \omega$ or $\NM(H)=\omega^{\alpha}$ for some ordinal $\alpha$, then $H$ is abelian-by-finite.
\end{theorem}

These results almost show \cite[Conjecture 5.31]{kr} which predicts a small, $nm$-stable compact $G$-group to be abelian-by-finite. The only remaining step is to prove Conjecture $(C')$ for infinite ordinal $\NM$-ranks.

This is related to another issue. In \cite[Conjecture 6.3]{kr}, the first author proposed the $\NM$-gap conjecture (formulated earlier by Newelski in the context of small profinite structures under the name of ${\cal M}$-gap conjecture \cite[Conjecture 2.2]{ne02}) which says that in a small Polish structure, the $\NM$-rank of every orbit is either a finite ordinal or $\infty$. In \cite{wa03}, this conjecture was proved for small, $nm$-stable profinite groups, i.e. each small, $nm$-stable profinite group is of finite $\NM$-rank. If this was true for small, $nm$-stable compact $G$-groups, then \cite[Conjecture 5.31]{kr} would be proved by Theorem \ref{Theorem XX}. However, in Section 3, we give a counter-example to the $\NM$-gap conjecture for small, $nm$-stable compact $G$-groups, i.e.\ an example of a small, $nm$-stable compact $G$-group of infinite ordinal $\NM$-rank. More precisely, we prove

\begin{theorem}
For every $\alpha <\omega_1$ there exists a small, $nm$-stable compact $G$-group $(H,G)$ such that $\NM(H)=\alpha$ and $H$ is abelian.
\end{theorem}

\section{Definitions and facts}

In the whole paper, the abbreviations $c$, $o$, $m$, $nm$ and $nwd$ come from the topological terms `closed', `open', `meager', `non-meager' and `nowhere dense', respectively.

We recall various definitions and facts from \cite{kr}. Assume that $(X,G)$ is a Polish structure. For $Y\subseteq X^n$ we define $Stab(Y):=\{ g \in G: g[Y]=Y\}$. We say that $Y$ is {\em invariant} [over a finite set $A$] if $Stab(Y)=G$ [$Stab(Y) \supseteq G_A$, respectively].

The {\em imaginary extension}, denoted by $X^{eq}$, is the union of all sets of the form $X^n/E$ with $E$ ranging over all invariant equivalence relations such that for all $a \in X^n$, $Stab([a]_E) \le_{c} G$. The sets $X^n/E$ are called the {\em sorts} of $X^{eq}$.

A subset $D$ of $X$ (or, more generally, of any sort of $X^{eq}$) is said to be {\em definable over a finite subset $A$} of $X^{eq}$ if $D$ is invariant over $A$ and $Stab(D) \le_c G$. We say that $D$ is {\em definable} if it is definable over some $A$. We say that $d \in X^{eq}$ is a {\em name} for $D$ if for every $f \in G$ we have $f[D]=D \iff f(d)=d$. 

It turns out that each set definable in $(X,G)$ [or in $X^{eq}$] has a name in $X^{eq}$. Moreover, working in $X^{eq}$, $\nmind$ has all the properties of an independence relation: G-invariance, symmetry, transitivity, and assuming smallness, the existence of independent extensions. We should mention here one more useful property of $\nmind$. For a finite $A \subseteq X^{eq}$ we define  $Acl(A)$ [$Acl^{eq}(A)$, respectively] as the collection of all elements of $X$ [$X^{eq}$, respectively] with countable orbits over $A$. Then, $a \in Acl^{eq}(A)$ iff $a \nmind_A B$ for every finite $B\subseteq X^{eq}$. In particular, $\NM(a/A)=0$ iff $a \in Acl^{eq}(A)$. We also define $dcl(A)$ as the set of elements of $X^{eq}$ which are fixed by $G_A$.

\defn A {\em compact [topological] $G$-space} is a Polish structure $(X,G)$, where $G$ acts continuously on a compact [topological] space $X$.\edefn

If $(X,G)$ is a compact $G$-space, we say that $D\subseteq X^n$ is {\em $A$-closed} for a finite $A\subseteq X^{eq}$ if it is closed and invariant over $A$. We say that it is {\em $*$-closed}, if it is $A$-closed for some $A$. 

It would be good to know that if $(X,G)$ is a topological $G$-space, then various topological properties of $X$ transfer to all the sorts of $X^{eq}$, but in general this is not true. That is the case, however, while working in a compact $G$-space $(X,G)$, provided that we restrict the collection of sorts in an appropriate way. This was a motivation for defining $X^{teq}$ ({\em topological imaginary extension}) as the disjoint union of the spaces $X^n/E$ with $E$ ranging over all $\emptyset$-closed equivalence relations on $X^n$. The spaces $X^n/E$ are called {\em topological sorts} of $X^{teq}$. Notice that each topological sort of $X^{teq}$ is, of course, a sort of $X^{eq}$.  It turns out that each topological sort $X/E$ together with the group $G/G_{X/E}$ is a compact $G$-space. If $E$ is $A$-closed for some finite set $A$, then after  replacing $G$ by $G_A$, $X/E$ can also be treated as a topological sort.

Now, we recall some facts about groups. We restrict ourselves to the situation
relevant to this paper, i.e. to small compact $G$-groups.

Let $(X,G)$ be a small compact $G$-space. We say that a group $H$ is {\em $*$-closed} in $X^{teq}$ if both $H$ and the group operation on $H$ are $*$-closed in $X^{teq}$. 
We would like to draw the reader's attention to the fact that we will use the same notation for several different operations: for $h_1,h_2 \in G$, $h_1h_2$ is the product computed in $H$; for $g_1,g_2 \in G$, $g_1g_2$ is the product computed in $G$; for $g \in G$ and $h \in X^{eq}$, $gh$ is the image of $h$ under $g$.

Assume $H$ is $*$-closed in $X^{teq}$; for simplicity, $\emptyset$-closed (then $(H,G/G_H)$ is a small compact $G/G_H$-group). Let $a \in H$ and $A\subseteq X$ be finite. We say that the orbit $o(a/A)$ is {\em $nm$-generic} (or that $a$ is {\em $nm$-generic} over $A$) if for all $b \in H$ with $a\nmind_A b$ one has that $b a\nmind A,b$. It turns out that $nm$-generics satisfy all the basic properties that generics satisfy in groups definable in simple theories, including existence. More precisely, \cite[Proposition 5.5]{kr} tells us that $o(a/A)$ is $nm$-generic iff $o(a/A) \subseteq_{nm} H$, and by smallness, such an orbit exists. 
Recall that the $\NM$-{\em rank} of an $A$-closed subset $D$ of $X^{teq}$ is the supremum of $\NM(d/A)$, $d \in D$. 
In the case when $\NM(H) <\infty$, the orbit $o(a/A)$ is $nm$-generic iff it is of maximal $\NM$-rank (equal to $\NM(H)$) among the orbits of elements of $H$. Using all of this, we get the following Lascar inequalities for groups and their corollary (see \cite[Corollary 5.6]{kr}).

\begin{fact}[Lascar inequalities for groups]
Let $(X,G)$ be a small compact $G$-space and $H_1\le H$ groups which are $*$-closed in $X^{teq}$. Then, $$\NM(H_1)+\NM(H/H_1) \leq \NM(H) \leq \NM(H_1) \oplus \NM(H/H_1).$$  
\end{fact} 

\begin{corollary}
Let $(X,G)$ be a small compact $G$-space and $H_1\le H_2$ groups which are $*$-closed in $X^{teq}$. Then:\\
(i) if $H_1 <_{nwd} H_2$ and $\NM(H_2) < \infty$, then $\NM(H_1)< \NM(H_2)$,\\
(ii) if $H_1 \le_{o} H_2$, then $\NM(H_1)= \NM(H_2)$.
\end{corollary}

Typical examples of Polish compact $G$-spaces are pairs $(X,Homeo(X))$, where $X$ is a compact metric space and $Homeo(X)$ is its group of homeomorphisms equipped with the compact-open topology. Analogously, typical examples of compact $G$-groups are pairs $(H,Aut(H))$, where $H$ is a compact metric group and $Aut(H)$ is its group of topological automorphisms equipped with the compact-open topology.

It was noticed in \cite[Proposition 5.7]{kr} that if $(H,G)$ is a small compact $G$-group, then $H$ is locally finite. Recall \cite[Proposition 5.12]{kr}.   

\begin{fact}
Let $(H,G)$ be a small compact $G$-group.\\
(i) $H$ is profinite, i.e. the inverse limit of a countable system of finite groups. 
Thus, it has a countable basis of open neighborhoods of identity consisting of clopen subgroups, and it has a countable basis consisting of clopen sets.\\
(ii) If $G=Aut(H)$ is the group of all topological automorphisms of $H$, then the topology on $G$ must coincide with the compact-open topology.\\
(iii) If $G$ is equipped with the compact-open topology, then $G$ has a countable basis of open neighborhoods of $id$ consisting of open subgroups. Thus, $G$ is isomorphic to a closed subgroup of $S_{\infty}$.
\end{fact}

The next fact consists of \cite[Remarks 5.10 and 5.11]{kr}.

\begin{fact}\label{clopen subsets}
Assume $(H,G)$ is a small compact $G$-group.\\
(i) For every clopen subset $U$ of $H$, $Stab(U)$ is a clopen subgroup of $G$, and so $[G:Stab(U)] \leq \omega$. In particular, this applies in the case when $U$ is an open subgroup of $H$.\\
(ii) If $G_0$ is a closed subgroup of $G$ of countable index, then $(H,G_0)$ is small, and $\nmind$ as well as $\NM$-rank computed in $(H,G)$ are the same as in $(H,G_0)$.
\end{fact}

Using exactly the same proofs, we get the following generalization of this fact.

\begin{remark}\label{clopen subsets 1}
(i) Let $(X,G)$ be a small compact $G$-space. Suppose $Y$ is a $\emptyset$-closed subset of $X^{teq}$. Then, for every clopen subset $U$ of $Y$, $Stab(U)$ is a clopen subgroup of $G$, and so $[G:Stab(U)] \leq \omega$.\\
(ii) Let $(X,G)$ be a small Polish structure. If $G_0$ is a closed subgroup of $G$ of countable index, then $(X,G_0)$ is small, and $\nmind$ as well as $\NM$-rank computed in $(X,G)$ are the same as in $(X,G_0)$.
\end{remark}

The following remark follows from \cite[Remark 3.10]{kr}.

\begin{remark}\label{definable subsets}
Let $(X,G)$ be a small compact $G$-space and $D$ an $A$-closed subset of $X^{teq}$. Then $(D,G_A/G_{AD})$ is a small compact $G$-space, and for tuples and subsets of $D$, $nm$-independence computed in $(X,G_A)$ coincides with $nm$-independence computed in $(D,G_A/G_{AD})$.
\end{remark}

In the $nm$-stable case, additionally we have \cite[Corollary 3.19]{kr}.

\begin{fact}\label{preservation of NM-rank}
Assume a small compact $G$-space $(X,G)$ is $nm$-stable. Let $D$ be an $A$-closed subset of $X^{teq}$, and $a,B \subseteq D$ finite. Then $\NM(a/AB)$ computed in $(X,G)$ is equal to $\NM(a/B)$ computed in $(D,G_A/G_{AD})$.
\end{fact}   

For this paper, the following two results, mentioned already in the introduction, are particularly important (see \cite[Theorem 5.19]{kr} and \cite[Theorem 5.24]{kr}).

\begin{fact}\label{Con A'}
If $(H,G)$ is a small, $nm$-stable, compact $G$-group, then $H$ is solvable-by-finite.
\end{fact}

\begin{fact}\label{Con B'}
If $(H,G)$ is a small compact $G$-group of finite $\NM$-rank,  then $H$ is nilpotent-by-finite.
\end{fact}

Let us finish this section with an easy remark which will be useful in Section 2.
\begin{remark}\label{remark 2.1}
Let $H$ be a topological space and $G$ a group acting on $H$ by homeomorphisms. Let $a \in H$ be such that the orbit $Ga$ is non-meager in $H$. Suppose $G_0$ is a countable index subgroup of $G$. Then the orbit $G_0a$ is also non-meager in $H$. 
\end{remark}
{\em Proof.} There are $g_i$, $i \in \omega$, in $G$ such that $G=\bigcup_{i \in \omega} g_iG_0$. So by assumption, we see that $\bigcup_{i \in \omega} g_iG_0a \subseteq_{nm} H$. Hence, $g_iG_0a \subseteq_{nm} H$ for some $i \in \omega$. This implies that $G_0a \subseteq_{nm} H$. \hfill $\blacksquare$

\section{Structure of small, $\boldsymbol{nm}$-stable compact $\boldsymbol{G}$-groups} 

In the first part of this section, we prove Conjecture $(B')$ and in consequence Theorem \ref{Theorem X} (repeated here as Theorem \ref{Thm 1}). In the second part, we show Conjecture $(C')$ (assuming additionally that $\NM(H) \leq \omega$) and in consequence Theorem \ref{Theorem XX} (here as Theorem \ref{Thm 2}), and we also obtain some other results.

The proof of \cite[Theorem 5.24]{kr} (i.e.\ Fact \ref{Con B'}) consists of \cite[Lemma 5.25]{kr}, \cite[Lemma 5.26]{kr} and the final conclusion. The assumption that $\NM$-rank is finite was necessary in the proof of \cite[Lemma 5.26]{kr} and in the final conclusion, because we did not know whether an increasing sequence of $A$-closed subgroups must stabilize (such a chain condition is present in small profinite groups, and that is why the proof in \cite{ne01} works for an arbitrary ordinal $\NM$-rank). Theorem \ref{counter-example} shows that such a chain condition fails in the class of small, $nm$-stable compact $G$-groups. Thus, we need another argument to show Conjecture $(B')$ in its full generality.

Below we will prove a strengthening of \cite[Lemma 5.25]{kr}, which will allow us to prove Conjecture $(B')$ without using any chain conditions. 

Recall the definition of an iterated centralizer: If a group $G$ acts by automorphisms on a group $H$, we put $C_{0,H}(G)=\{e\}$, and $$C_{n+1,H}(G)=\{h\in N_H(C_{n,H}(G)):Gh\subseteq hC_{n,H}(G)\}.$$ 
Then $C_{1,H}(G)=C_H(G)$ is the subgroup of points in $H$ fixed by $G$, and $C_{n+1,H}(G)$ the subgroup of points normalizing $C_{n,H}(G)$, and fixed by $G$ modulo $C_{n,H}(G)$.

\begin{lemma}\label{abelian groups}
Suppose $D,H$ are abelian groups $*$-closed in $X^{teq}$ for some small compact $G$-space $(X,G)$, and $\NM(H) < \infty$. Assume that $D$ acts by automorphisms on $H$, and the action is $*$-closed. Then there exists a clopen subgroup $F$ of $D$ such that $C_{n,H}(F)$ is clopen in $H$ for some $n \in \omega$.
\end{lemma}
{\em Proof.} Wlog everything is invariant over $\emptyset$.  
When we say that an element is $nm$-generic (without specifying parameters), we mean that it is $nm$-generic over $\emptyset$. Note that $C_{n,H}(F)$ will be invariant under $D$, as $D$ is abelian.
First we need two claims.\\[3mm]
{\bf Claim 1.} If there is an $nm$-generic $a$ in $H$ for which $Da$ is finite, then some clopen subgroup $D_0$ of $D$ fixes pointwise a clopen subgroup of $H$.\\[2mm]
{\em Proof of Claim 1.} By assumption, there is a clopen subgroup $D_0$ of $D$ such that $a \in C_H(D_0)$. By Remark \ref{clopen subsets 1}, $Stab_G(D_0)$ is a subgroup of countable index in $G$. 
In virtue of the assumption that $a$ is $nm$-generic, $o(a) \subseteq_{nm} H$, and so by Remark \ref{remark 2.1},
$$o_1(a):=Stab_G(D_0)a \subseteq_{nm} H.$$
We also have $o_1(a) \subseteq C_H(D_0)\le_c H$. So, we get that $C_H(D_0)$ is a clopen subgroup of $H$. \hfill $\square$\\

The proof of the second claim is a modification of the proof of \cite[Lemma 5.25]{kr}. \\[3mm]
{\bf Claim 2}
Consider any $a \in H$ with $Da$ infinite. Let $\alpha$ be maximal such that  $\NM(\overline{\langle Da \rangle}) \geq \omega^\alpha$. Let $H_1$ be a $*$-closed subgroup of $\overline{\langle Da \rangle}$ invariant under $D$ of minimal $\NM$-rank greater or equal to $\omega^\alpha$. Then some open subgroup of $D$ fixes pointwise a clopen subgroup of $H_1$.\\[3mm] 
{\em Proof of Claim 2.} 
By the choice of $H_1$, for every $*$-closed subgroup $H_2$ of $H_1$ invariant under $D$ either $[H_1:H_2]<\omega$ or $\NM(H_2) <\omega^\alpha$.

Let $R$ and $S$ be the rings of endomorphisms of $H_1$ and $\overline{\langle Da \rangle}$, respectively, generated by $D$. Then $R$ and $S$ are commutative, and $R \cong S/J$ for some $J \lhd S$. As in the proof of \cite[Lemma 5.25]{kr}, S is locally finite, and hence $R$ is locally finite, too.

Consider any $r \in R$. Then $r[H_1]$ and $ker(r)$ are $*$-closed subgroups of $H_1$ invariant under $D$.
Thus, by the choice of $H_1$, either $[H_1:r[H_1]]<\omega$ or $\NM(r[H_1])<\omega^\alpha$. The second case implies that $\NM(ker(r)) \geq \omega^\alpha$, so $[H_1: ker(r)]<\omega$, and hence $r[H_1]$ is finite.

As in the proof of \cite[Lemma 5.25]{kr}, the above dichotomy implies that $I:=\{ r \in R: r[H_1]\; \mbox{is finite}\}$ is a prime ideal of $R$, and so $R/I$ is a locally finite integral domain. Thus, $R/I$ is countable.
 
Put $\widetilde{D}:=D/D_{H_1}\le R^*$. Then, there are $d_i \in \widetilde{D}$, $i \in \omega$, such that $\widetilde{D} \subseteq \bigcup_i d_i+I$. 

Let $\{A_j: j \in \omega\}$ be a basis of open neighborhoods of $e$ in $H_1$ consisting of clopen subgroups. Then, for every $d \in \widetilde{D}$ there is $i \in \omega$ and $j \in \omega$ such that $A_j \subseteq ker(d-d_i)$. Hence, 
$$\widetilde{D}=\bigcup_i \bigcup_j \{ d \in \widetilde{D}: (d-d_i) \upharpoonright A_j = 0\}.$$
By the Baire category theorem applied to the compact group $\widetilde{D}$, one can find $i,j \in \omega$ such that the set $\{ d \in \widetilde{D}: (d-d_i) \upharpoonright A_j = 0\}$ is non-meager in $\widetilde{D}$. But this set is also closed, so it has a nonempty interior $U$.
Then, for every $b \in A_j$, $Ub=\{ d_ib\}$, which implies $d_i^{-1}U \subseteq \widetilde{D}_b$. Therefore, $d_i^{-1}U \subseteq \widetilde{D}_{A_j}$, and so $\widetilde{D}_{A_j}$ is clopen in $\widetilde{D}$. Putting $F:=\pi^{-1}[\widetilde{D}_{A_j}]$, where $\pi: D \to \widetilde{D}$ is the quotient map, we see that $F$ is a clopen subgroup of $D$ fixing $A_j$ pointwise. \hfill $\square$\\

Now, using Claims 1 and 2, we will complete the proof of Lemma \ref{abelian groups}.
Choose $a$ $nm$-generic in $H$. Put $L=\overline{\langle Da \rangle}$, an $a$-closed subgroup of $H$. If $L$ is finite, we are done by Claim 1. So, assume that $L$ is infinite.

Using Claim 2, we will construct sequences of closed groups
$$D=D_0\ge D_1 \ge \dots$$ and 
$$\{ e \}=L_0\le L_1\le \dots\le L$$ 
such that all $D_i$ are clopen in $D$, each $L_i$ is invariant under $D$, and
\begin{enumerate}
\item[(i)] if $\alpha_i$ is maximal with $\NM_i(\overline{\langle D(a +L_i) \rangle}) \geq \omega^{\alpha_i}$, then $\NM_{i+1}(C_{L/L_i}(D_{i+1})) \geq \omega^{\alpha_i}$ ($\NM_i$-rank is $\NM$-rank computed in $(X, Stab_{G_a}(D_0) \cap \dots \cap Stab_{G_a}(D_i))$),
\item[(ii)] $L_{i+1}=\sigma_i^{-1}[C_{L/L_i}(D_{i+1})]$, where $\sigma_i:L \to L/L_i$ is the quotient map.
\end{enumerate}

Suppose $D=D_0\ge D_1\ge \dots \ge D_i$ and $\{ e \}=L_0\le L_1\le \dots\le L_i$ have been constructed. If $L/L_i$ is finite, we put $D_{i+1}=D_i$ and $L_{i+1}=\sigma_i^{-1}[C_{L/L_i}(D_{i+1})]$. If $L/L_i$ is infinite, then there is a unique maximal $\alpha_i$ with $\NM_i(\overline{\langle D(a+L_i) \rangle}) \geq \omega^{\alpha_i}$. By Claim 2 (applied to $H/L_i$ instead of $H$) and Remark \ref{clopen subsets 1}, there is a clopen subgroup $D_{i+1}$ of $D$, which can be chosen to be a subgroup of $D_i$ (by taking the intersection with $D_i$), such that we have $\NM_{i+1}(C_{L/L_i}(D_{i+1})) \geq \omega^{\alpha_i}$. Define $L_{i+1}=\sigma_i^{-1}[C_{L/L_i}(D_{i+1})]$, and note that by commutativity of $D$, it is invariant under $D$. The construction is completed. 

Of course, $\overline{\langle D(a +L_i) \rangle}=L/L_i$.
By the Lascar inequalities for groups and Remark \ref{clopen subsets 1}, for every $i$ such that $L/L_i$ is infinite $\NM_{i+1}(L/L_{i+1})<\NM_{i}(L/L_i)$. So, after finitely many, say $n$, steps, $[L:L_n]<\omega$. Hence, $D_n(a+L_n) \subseteq L/L_n$ is finite.

Define recursively a sequence $\{ e \}=K_0\le K_1\le \dots\le K_n$ of closed subgroups of $H$ invariant under $D$ in the following way.\begin{enumerate}
\item[] $K_{i+1}:=\tau_i^{-1}[C_{H/K_i}(D_{i+1})]$, where $\tau_i:H \to H/K_i$ is the quotient map.
\end{enumerate}

Then, $L \cap K_i=L_i$ for every $i$, and we get that $D_n(a+K_n)$ is a finite subset of $H/K_n$.

Let $G^*=Stab_G(D_0) \cap \dots \cap Stab_G(D_n)$. Since $a$ is $nm$-generic in $H$,
by Remark \ref{remark 2.1}, $a$ is $nm$-generic in $(H,G^*/G^*_H)$. Since $K_n$ is invariant under $G^*$, this implies that $a+K_n$ is $nm$-generic in $(H/K_n, G^*/G^*_{H/K_n})$. 

Thus, by Claim 1 (applied to $H/K_n$ instead of $H$), there exists a clopen $D_{\omega}\le D_n$ such that $C_{H/K_n}(D_{\omega})$ is clopen in $H/K_n$. Put $K_{n+1}=\tau_n^{-1}[C_{H/K_n}(D_{\omega})]$. Then, $K_{n+1}$ is clopen in $H$. 

Define $F=D_{\omega}$. We easily check that $K_i\le C_{i,H}(F)$ for all $i \leq n+1$. Since $K_{n+1}$ is clopen in $H$, so is $C_{n+1,H}(F)$. \hfill $\blacksquare$\\

Using Fact \ref{Con A'} and Lemma \ref{abelian groups}, we will prove Conjecture $(B')$. 
In the proof, we use an easy observation that if $K$ is a [$A$-invariant] solvable subgroup of a compact $G$-group $(H,G)$, then the group $\overline{K}$ is [$A$-]closed and solvable of the same solvability class as $K$. A similar remark is true for nilpotent subgroups, too.

\begin{theorem}\label{Thm 1}
If $(H,G)$ is a small, $nm$-stable compact $G$-group, then $H$ is nilpotent-by-finite.
\end{theorem}
{\em Proof.} In virtue of Fact \ref{Con A'}, $H$ is solvable-by-finite. So, replacing $H$ by a clopen subgroup and using Fact \ref{clopen subsets} and Remark \ref{definable subsets}, we can assume that $H$ is solvable. We proceed  by induction on the solvability class of $H$. 

In the abelian case, we have the desired conclusion. Assume the conclusion holds for groups of solvability class smaller than the solvability class of $H$.
As $\overline{H'}$ is $\emptyset$-closed and of smaller solvability class than $H$, we get that $\overline{H'}$ is nilpotent-by-finite.
Replacing $H$ by a clopen subgroup and using Fact \ref{clopen subsets} and Remark \ref{definable subsets}, we can assume that $\overline{H'}$ is nilpotent.
Now, we do the induction on the nilpotency class of $\overline{H'}$.

We apply the same argument for the base and for the induction step. Assume the conclusion holds for groups whose commutator group is of smaller nilpotency class than $\overline{H'}$. We see that $H/\overline{H'}$ acting by conjugation on $Z(\overline{H'})$ satisfies the assumption of Lemma \ref{abelian groups}. So, there exists a clopen subgroup $F/\overline{H'}$ of $H/\overline{H'}$ for which there is $n$ with $[Z(\overline{H'}): C_{n,Z(\overline{H'})}(F/\overline{H'})]<\omega$. Then, $[Z(\overline{H'}): Z_n(F) \cap Z(\overline{H'})]<\omega$, 
where $Z_n(F)$ is the $n$-th term in the upper central series of $F$. 
Thus, we can choose a clopen subgroup $H_0$ of $F$ so that $Z(\overline{H'}) \cap H_0 = Z_n(F) \cap Z(\overline{H'}) \cap H_0$. Then, $H_0$ is a clopen subgroup of $H$. Moreover, since $Z_n(H_0) \geq Z_n(F) \cap H_0$, we see that $Z_n(H_0) \geq  Z(\overline{H'}) \cap H_0$. Using this and the fact that $H_0'/(Z(\overline{H'}) \cap H_0')$ is of smaller nilpotency class than $\overline{H'}$, we conclude that   $\overline{(H_0/Z_n(H_0))'}$ is of smaller nilpotency class than $\overline{H'}$ (it is trivial in the base induction step, i.e.\ if $\overline{H'}$ was abelian). Hence, after replacing $H$ by $H_0$ (which we can do by  Fact \ref{clopen subsets} and Remark \ref{definable subsets}), $H/Z_n(H)$ becomes nilpotent-by-finite by the induction hypothesis. Therefore, $H$ is nilpotent-by-finite. \hfill $\blacksquare$

\begin{proposition}\label{finite exponent}
If $(H,G)$ is a small compact $G$-group and $H$ is solvable-by-finite, then $H$ has finite exponent. In particular, if $(H,G)$ is a small, $nm$-stable compact $G$-group, then $H$ has finite exponent.
\end{proposition}
{\em Proof.} If $(H,G)$ is a small, $nm$-stable compact $G$-group, then in virtue of Theorem \ref{Thm 1}, $H$ is nilpotent-by-finite. So, it is enough to prove the first statement of the proposition. Since $H$ is solvable-by-finite, it is enough to prove the assertion assuming that $H$ is abelian.

Choose $h \in H$ $nm$-generic. Since $H$ is locally finite \cite[Proposition 5.7]{kr}, $h$ has finite order $n$. 
As $H[n]:=\{ a \in H: na=0 \}$ contains $o(h)$, $H[n]$ is a non-meager, and hence 
clopen subgroup of $H$. Thus, $H$ has finite exponent. \hfill $\blacksquare$

\begin{proposition}\label{big centralizers}
Let $(H,G)$ be a small compact $G$-group. Assume that there is an $nm$-generic (over $\emptyset$) $h$ with open centralizer. Then $H$ is abelian-by-finite.
\end{proposition}
{\em Proof.}  
Using our assumption together with  Remarks \ref{clopen subsets 1} and \ref{remark 2.1}, we easily get
\begin{equation*}\tag{$*$}
o_1(h):=Stab_G(C(h))h \subseteq_{nm} H.
\end{equation*}

Since $h \in C(h)$, we have $o_1(h) \subseteq C(h)$. In fact, as $h$ centralizes $C(h)$, all of  $o_1(h)$ centralizes $C(h)$. Thus, $o_1(h) \subseteq Z(C(h))$. Hence, $Z(C(h))$ is open in $H$ by $(*)$. Of course, it is also abelian. \hfill $\blacksquare$\\

Our second goal is to prove Conjecture $(C')$ for groups of $\NM$-rank at most $\omega$. The proof is a modification and simplification of the proof from \cite{wa03} (a simplification because we omit the application of Schlichting's theorem). We will usually skip the parts of the proofs which are identical to \cite{wa03}.

Recall that two subgroups $H_1$ and $H_2$ of a given group $H$ are said to be {\em commensurable}, denoted 
$H_1 \sim H_2$, if their intersection $H_1 \cap H_2$ has finite index both in $H_1$ and in $H_2$. 

Now, we recall the notion of a minimal group and virtual isogeny.
Let $(X,G)$ be a small compact $G$-space and $H$ a group $*$-closed in $X^{teq}$. 

\defn
A subgroup $A$ of $H$ is {\em minimal} if it is infinite, $*$-closed, and every $*$-closed subgroup of infinite index in $A$ is finite.
\edefn

\defn
Let $A$ and $B$ be abelian, minimal subgroups of $H$. A {\em virtual isogeny} $f$ between $A$ and $B$ is a $*$-closed isomorphism $f:D/K \to I/C$, where $D$ is open in $A$, $I$ is open in $B$, and $K$, $C$ are both finite. Two virtual isogenies  $f_1$ and $f_2$ are {\em equivalent}, denoted $f_1 \sim f_2$, if the graphs of $f_1$ and $f_2$ (treated as subsets of $A \times B$) are commensurable (equivalently, if the derived maps from $D_1 \cap D_2$ to $(I_1 + I_2)/(C_1+C_2)$ agree on an open subgroup of $A$).
\edefn

It is standard that in an abelian, minimal subgroup $A$ of $H$, the family of virtual autogenies (isogenies from $A$ to $A$), with addition and composition as operations, forms the set of non-zero elements of a division ring $R$.

In the sequel, we will use the following notation: if $D_a$ is an $a$-closed set and $b=ga$ for some $g \in G$, then $D_b:=g[D_a]$. 

The following Lemma is a variant of \cite[Lemma 14]{wa03}. The proofs differ because in our context, the algebraic closure, Acl, of a finite set does not need to be finite.

\begin{lemma}\label{R is locally finite}
Let $(X,G)$ be a small compact $G$-space and $H$ a group $*$-closed in $X^{teq}$. Let $A$ be an abelian, minimal subgroup of $H$. Let $R$ be the division ring considered above, i.e. the ring whose non-zero elements are virtual autogenies of $A$.
Then $R$ is locally finite. Moreover, if $P$ is a finite set of parameters such that $H$ and $A$ are $P$-closed, then for any finite tuple  $a$ of parameters from $X$ and for every $a$-closed virtual autogeny $f_a$ of $A$ the equivalence relation $E(x,y)$ on $o(a/P)$ given by $f_x \sim f_y$ has finitely many classes. 
\end{lemma} 
{\em Proof.}
Assume for simplicity that $P=\emptyset$.
In virtue of Proposition \ref{finite exponent}, $A$ has finite exponent. Thus, by minimality of $A$, there is a prime number $p$ such that
$A[p]:=\{c \in A: pc=0\}$ is $\emptyset$-closed and open in $A$. Replacing $A$ by $A[p]$, we do not change $R$. So, we can assume that $A$ is elementary abelian of exponent $p$.

Let $\{f_1,\dots, f_n\}$ be a finite set of virtual autogenies of $A$. 
Let $C_1, \dots, C_n$ be finite subgroups of $A$ such that $f_i \cap (\{ 0 \} \times A)= \{ 0 \} \times C_i$ (where $f_i$ are treated as subsets of $A \times A$). Since $A$ is profinite,
we can choose a clopen $N\lhd A$ with $N \cap C_i =\{ e \}$ for every $i$. Put $f_i' = f_i \cap (A \times N)$. Then, each $f_i'$ is a homomorphism from some clopen $A_i\le A$ to $A$, and $f_i'\sim f_i$ (notice that since $N$ does not need to be $*$-closed, $f_i'$'s are not necessarily virtual autogenies; nevertheless, $f_i \sim f_i'$ has a perfect sense).

Since $A$ is elementary abelian, each $A_i$ has a finite complement $B_i$ in $A$. Put $f_i''= f_i'\oplus (B_i \times \{ 0\})$. Then, each $f_i''$ is an endomorphism of $A$ with finite kernel, and $f_i'' \sim f_i$.

Let $\widetilde{R}$ be the ring of endomorphisms of $A$, and $\langle \overline{f''} \rangle$ the subring of $\widetilde{R}$ generated by  $\overline{f''}:=(f_1'',\dots,f_n'')$. In this paragraph, we work in the small Polish structure $(X, Stab_G(N))$. Then, all $f_i''$ are $*$-closed. Let $\overline{a''}$ be a finite tuple over which $\overline{f''}$ is defined. Choose $h\in A$ $nm$-generic over $\overline{a''}$. Then, for any $f,f' \in \langle \overline{f''} \rangle$ we have that $f(h)$ and $f'(h)$ belong to $dcl(\overline{a''},h) \cap A$. Since $dcl(\overline{a''},h) \cap A$ is a closed subgroup of $A$,  it must be finite by smallness of $(X, Stab_G(N))$.
On the other hand, if $f(h)=f'(h)$, then $h \in ker(f-f')$, so using the fact that $ker(f-f')$ is $\overline{a''}$-closed and $h$ is $nm$-generic over $\overline{a''}$, we get $f \sim f'$. We conclude that modulo $\sim$, $\langle \overline{f''} \rangle$ is finite. Thus, the subring of $R$ generated by $\{[f_1]_{\sim},\dots,[f_n]_\sim\}$ is also finite. Hence, $R$ is locally finite, and so it is a locally finite field. Since all conjugates of $[f_a]_\sim$ are of the same finite order in this field, $E$ has finitely many classes. \hfill $\blacksquare$\\

Now, we will prove a variant of \cite[Theorem 15]{wa03}, which is strong enough to prove Conjecture $(C')$ for groups of $\NM$-rank at most $\omega$. It has a weaker conclusion than \cite[Theorem 15]{wa03}, and this allows us to eliminate the application of Schlichting's theorem from the proof.

\begin{proposition}\label{commensurability} 
Let $(X,G)$ be a small compact $G$-space and $H$ a group $\emptyset$-closed in $X^{teq}$. Assume $H$ is abelian and $\NM(H) < \omega$. Let $H_1$ be an $h_1$-closed subgroup of $H$ for some finite tuple $h_1$ of parameters from $X$. Then there is a a finite tuple $k$ of parameters from $X$ with $k \nmind h_1$ and a $k$-closed subgroup $K$ of $H$ which is commensurable with $H_1$.
\end{proposition} 
{\em Proof.} As in \cite{wa03}, first we prove the statement for minimal groups.\\[3mm]
{\bf Claim} Let $A$ be a minimal, $*$-closed subgroup of $H$; say it is $a$-closed for some finite tuple $a$ of parameters from $X$ (so we can write $A=A_a$).  Then there is $b \in o(a)$ such that $b \nmind a$ and $A_a \sim A_b$.\\[2mm]
{\em Proof of Claim.} The first two paragraphs of the proof are the same as in \cite{wa03}, so we only give the conclusion coming from them: there is a finite tuple $\overline{a}$ of parameters, a conjugate $A_i$ of $A$ which is ${\overline{a}}$-closed, $a' \in o(a)$ with $a' \nmind \overline{a}$, and finitely many $\{\overline{a}, a'\}$-closed virtual autogenies $f_{a',j}$, $j=1,\dots,k$, of $A_i$ such that for any $x, y \in S:=o(a'/\overline{a})$, if $f_{x,j} \sim f_{y,j}$ for all $j$, then $A_{x} \sim A_{y}$. 

For any $x \in S$ choose $g_x \in G$ such that $g_x a=x$.
Let $(x_i)_{i \in \omega}$ be a sequence of elements of $S$ such that $x_i \nmind_{\overline{a}} x_{<i}$. 

Let $F_j(x,y)$ be the equivalence relation on $S$ given by $f_{x,j} \sim f_{y,j}$.
By Lemma \ref{R is locally finite}, each $F_j$ has finitely many classes. So, there exist $i_1<i_2$ such that $A_{x_{i_1}} \sim A_{x_{i_2}}$. Put $b=g_{x_{i_1}}^{-1}(x_{i_2})$. Then, we have: $b \in o(a)$, $b \nmind a$ and $A_a \sim A_b$. \hfill $\square$\\

Now, we will prove the proposition by induction on $\NM(H_1)$. If $H_1$ is minimal, we are done by the Claim. For the induction step choose $A_a\le H_1$ which is minimal (and $a$-closed). Then, $A_a$ is also $\{ a,h_1\}$-closed, so
using the Claim, we can find $b \nmind a,h_1$ such that $A_b \sim A_a$. Put $H_1^*=H_1A_b/A_b\cong H_1/(H_1 \cap A_b)$. We see that $H_1^*$ is $h_1$-closed (in $X^{teq}$ computed in $(X,G_b)$) and of smaller $\NM$-rank than $H_1$. So, by the induction hypothesis, 
there is a finite tuple $k_1$ of parameters from $X$ with $k_1 \nmind_{b} h_1$ and a $k_1$-closed subgroup $K_1$ of $H/A_b$ (still working in $(X,G_b)$) such that $K_1 \sim H_1^*$. Let $\pi: H \to H/A_b$ be the quotient map, and $K:=\pi^{-1}[K_1]$. Then, $K$ is $\{k_1,b\}$-closed and $K \sim H_1$. 
Since $b \nmind h_1$ and $k_1 \nmind_{b} h_1$,  we also see that $k_1,b \nmind h_1$. So, the tuple $k:=(k_1,b)$ satisfies the requirements. \hfill $\blacksquare$

\begin{theorem}\label{Thm 2}
If $(H,G)$ is a small compact $G$-group, and either $\NM(H) < \omega$ or $\NM(H)=\omega^{\alpha}$ for some ordinal $\alpha$, then $H$ is abelian-by-finite.
\end{theorem}
{\em Proof.} By Theorem \ref{Thm 1}, $H$ is nilpotent-by-finite. So, using Fact \ref{clopen subsets}, Remark \ref{definable subsets} and Fact \ref{preservation of NM-rank}, we can assume that $H$ is nilpotent. Next, by induction on the nilpotency class of $H$, one can easily reduce the situation to the case when $H$ is of nilpotency class 2.

First, we prove the theorem in the case when $\NM(H)<\omega$.
As in \cite{wa03}, for $g \in H$ put $H_g:= \{ (hZ(H),[h,g]): h \in H\}$, a $g$-closed subgroup of the abelian group $H/Z(H) \times Z(H)$. Choose $nm$-generic $g\in H$.
By Proposition \ref{commensurability}, $H_g \sim K_k$ 
for some $k$-closed subgroup $K_k$ with $k \nmind g$. Choose $g' \in o(g/k)$ with $g' \nmind_k g$. Then, $H_{g'} \sim K_k \sim H_g$, so 
\begin{equation*}\tag{$*$}
H_{g'} \sim H_g.
\end{equation*}
Moreover,
\begin{equation*}\tag{$**$}
g' \nmind g.
\end{equation*}

Let $H_1=\pi^{-1}[\pi_1[H_{g'}\cap H_g]]$, where $\pi: H \to H/Z(H)$ is the quotient map and $\pi_1: H/Z(H) \times Z(H) \to H/Z(H)$ is the projection on the first coordinate. For $h \in H_1$, $[h,g]=[h,g']$, whence $[h,g'g^{-1}]=e$, which means that $h \in C(g'g^{-1})$. Thus, 
$H_1 \le C(g'g^{-1})$.
On the other hand, $[H:H_1]<\omega$  by $(*)$. So, we conclude that $C(g'g^{-1})$ is a closed subgroup of finite index in $H$, whence open in $H$.  Moreover, by $(**)$, we have that $g'g^{-1}$ is $nm$-generic. Using these two observations together with Proposition \ref{big centralizers}, we get that $H$ is abelian-by-finite.

Now, consider the case $\NM(H)=\omega^\alpha$. If $\NM(Z(H)) = \omega^\alpha$, we are done. Otherwise, for any $g \in H$ define a homomorphism $f_g: H \to Z(H)$ by $f_g(h)=[h,g]$. Then, $\NM(Im(f_g)) \leq \NM(Z(H))<\omega^\alpha$, so $\NM(ker(f_g))=\omega^\alpha$. Thus, $ker(f_g)$ has finite index in $H$, which implies that $[H:C(g)]<\omega$. Since this holds for every $g \in H$ (in particular, for $nm$-generic $g$), we are done by Proposition \ref{big centralizers}. \hfill $\blacksquare$\\

Although the main goals of Section 2 have already been achieved, something interesting about small compact $G$-groups can still be said.

Proposition 2.8 has a weaker conclusion than \cite[Theorem 15]{wa03}, but strong enough to prove Theorem \ref{Thm 2}. Nevertheless, \cite[Theorem 15]{wa03} is interesting in its own right, and it would be nice to know whether it holds for small compact $G$-groups. Below we prove a variant of it. 
However, in the proof we will need Schlichting's theorem, which we recall now.

\begin{fact}\label{Schlichting}
Let $G$ be any group, and ${\EuFrak h}$ a family of uniformly commensurable subgroups. Then there is a subgroup $N$ of $G$, a finite extension of a finite intesection of groups in ${\EuFrak h}$ (and hence commensurable with them), such that $N$ is invariant under all automorphisms of $G$ fixing ${\EuFrak h}$ setwise.
\end{fact}

\begin{theorem}\label{Acl(emptyset)}
Let $(X,G)$ be a small compact $G$-space. Assume additionally that $G$ is equipped with the compact-open topology (or any topology having a basis at $id$ consisting of open subgroups). Let $H$ be a group $\emptyset$-closed in $X^{teq}$, and suppose that $H$ is abelian and $\NM(H)<\omega$. Then, for every $*$-closed $H_1\le H$ there is an $Acl^{eq}(\emptyset)$-closed subgroup $K$ of $H$ commensurable with $H_1$.\end{theorem}

In particular, the theorem applies in the situation when 
$G=Homeo(X)$, or $G=Aut(X)$ where $X$ is a compact metric group.\\[3mm]
{\em Proof.} By the inductive argument (on $\NM(H_1)$)  in the last paragraph of the proof of \cite[Theorem 15]{wa03}, we see that it is enough to show the assertion for $H_1$ minimal. To have the same notation as in \cite{wa03}, put $A=A_a:=H_1$. 

The Claim in the proof of Proposition \ref{commensurability} yields  $a' \in o(a)$ such that $a' \nmind a$ and $A_{a} \sim A_{a'}$.

In order to finish the proof, we need an extra topological argument as we do not have a generalization of \cite[lemma 11]{wa03}.

We define:
$$
\begin{array}{ll}
C=\{ b \in o(a): A_b \sim A\}, & \\
C_n^-=\{b \in o(a): [A: A\cap A_b] \leq n\}, & \\
C_n^+=\{ b \in o(a): [A_b:A\cap A_b]\leq n\},\\ 
C_n=C_n^- \cap C_n^+. &
\end{array}
$$
Then $C= \bigcup_{n} C_n$.\\[3mm]
{\bf Claim} $P_n:= \{ g \in G : ga \in C_n\}$ is closed.\\[3mm]
{\em Proof of Claim.} We have $P_n= P_n^- \cap P_n^+$,
where $P_n^- = \{ g \in G: ga \in C_n^-\}$ and  $P_n^+ = \{ g \in G: ga \in C_n^+\}$. Moreover, $P_n^+ = (P_n^-)^{-1}$. Thus, it is enough to show that $P_n^-$ is closed.

Suppose for a contradiction that $P_n^-$ is not closed, i.e. there are $g_k \in P_n^-$ (for $k \in \omega$) such that $g_k \longrightarrow g$, but $g \notin P_n^-$. Then, $ga \notin C_n^-$, which means that $[A:A \cap g[A]]>n$ (as $g[A]=A_{ga}$). Choose $D\subseteq A$ dense and countable.  Then, there are $d_0, \dots, d_n\in D$ such that $d_id_j^{-1} \notin g[A]$ for all $i \ne j$. Take an open neighborhood $U$ of $g[A]$ disjoint from $\{ d_id_j^{-1} : i,j=0,\dots, n\; \mbox{and}\; i \ne j \}$. Then, $V:=\{h \in G: h[A] \subseteq U\}$ is an open neighborhood of $g$. Thus, there is $k_0$ such that for all $k\geq k_0$, $g_k \in V$. Consider any $k \geq k_0$. Then, $g_k[A] \subseteq U$. But $g_k[A]=A_{g_{k}a}$. Whence, $d_id_j^{-1} \notin A_{g_{k}a}$ for all $i,j=0,\dots, n$ such that $i \ne j$. This implies that $[A:A \cap A_{g_{k}a}] >n$, a contradiction with the choice of $g_k$.\hfill $\square$\\

Since $A_a \sim A_{a'}$, we see that $\{g \in G: ga \in o(a/a')\} \subseteq P_n$ for some $n$. As $a \nmind a'$, directly from the definition of $\nmind$, we get that $\{g \in G: ga \in o(a/a')\}$ is non-meager. Thus, $P_n$ is non-meager.  So, by the Claim, $int(P_n) \ne \emptyset$.
Since $G$ has a basis of open neighborhoods of $id$ consisting of open subgroups (this is the only place where this extra assumption is used), there is a clopen subgroup $G_0$ of $G$ and $g \in G$ with $gG_0 \subseteq P_n$. This means that $gG_0a \subseteq C_n$, and so we can apply Fact \ref{Schlichting} to the uniformly commensurable family $\{ A_x: x \in gG_0a\}$ of subgroups of $H$. As a result, we obtain a $*$-closed subgroup $K$ of $H$ which is commensurable with $A$ and invariant under $Stab_{G}(gG_0a)$. But since the index of $G_0$ in $G$ is countable, 
$gG_0a$ has countably many different conjugates by the elements of $G$, i.e. $[G: Stab_{G}(gG_0a)]\leq \omega$. Thus, $K$ is $Acl^{eq}(\emptyset)$-closed. \hfill $\blacksquare$\\

It is worth mentioning that \cite[Theorem 15]{wa03} was used to prove the $\NM$-gap conjecture for small, $nm$-stable profinite groups. In our more general context of small, $nm$-stable compact $G$-groups, this is impossible because the conjecture is false as we will see in Section 3.

At the end of this section, we discuss one more issue. In the context of $\omega$-categorical, supersimple groups, the strongest statement says that they are finite-by-abelian-by-finite \cite{ew}. In the context of small profinite groups, or even small compact $G$-groups, being finite-by-abelian-by-finite immediately gives us abelian-by-finite by taking an appropriate clopen subgroup. (Indeed, if a profinite group $H$ has a finite index subgroup $H_1$ possessing a finite, normal subgroup $H_2$ such that $H_1/H_2$ is abelian, then choosing a clopen subgroup $H_0$ of $H$ so that $H_0 \cap H_2 = \{ e \}$, we get that $H_0 \cap H_1$ is abelian and of finite index in $H$.) However, one can ask what happens if we know that the group is countable-by-abelian-by-countable.

\begin{remark}
Let $(H,G)$ be a small compact $G$-group. If $H$ is countable-by-abelian-by-countable, then $H$ is abelian-by-finite.
\end{remark} 
{\em Proof.} By assumption, there are $H_1\le H$ and $H_2\lhd H_1$ such that $[H:H_1]\leq \omega$, $|H_2|\leq \omega$ and $H_1/H_2$ is abelian. By the Baire category theorem, we get $H_1 \subseteq_{nm} H$, whence by smallness, there is $h \in H_1$ $nm$-generic in $H$. Moreover, $H_1'\le H_2$, so $[H_1: H_1 \cap C(h)]\leq \omega$, which implies $[H:C(h)]\leq \omega$. Hence, $C(h)$ is open in $H$, and we finish using Proposition \ref{big centralizers}. \hfill $\blacksquare$

\section{Examples of groups of infinite ordinal $\boldsymbol{\NM}$-rank}

In this section, we give counter-examples to the $\NM$-gap conjecture, even in the class of small, $nm$-stable compact $G$-groups. The examples are pretty simple, but the proof that they satisfy the required properties is quite long and technical (as usual, when one needs to describe orbits in a concrete example).

\begin{remark}
If $(X,G)$ is a small Polish structure, then the set of values $\NM(o)$, with $o$ ranging over all orbits over finite sets, is the interval $[0,\alpha)$ (together with the element $\infty$, if $X$ is not $nm$-stable) for some countable ordinal $\alpha$.
\end{remark}
{\em Proof.} It is easy to see that the $\NM$-rank of an element of $X^{eq}$ over a finite subset of $X^{eq}$ is bounded by the $\NM$-rank of a finite tuple from $X$ over a finite subset of $X$. So, it is enough to work with finite tuples and subsets of the home sort.

It is clear that the set of values of $\NM$-rank  is an initial interval, possibly together with $\infty$. Since $\NM$-rank is invariant under $G$ and there are only countably many orbits on each $X^n$ 
(and so countably many orbits on the collection of all finite subsets of $X$), we see that $\NM$-rank takes countably many values.\hfill $\blacksquare$\\

The next theorem shows that the $\NM$-gap conjecture is false for small, $nm$-stable compact $G$-groups. Moreover, it shows that the ascending chain condition on $\emptyset$-closed subgroups fails. Recall, once again, that such a condition holds for small profinite groups \cite{ne01}, which is an important tool in proving structural theorems about small profinite groups. The lack of this chain condition in our context is one of the reasons why in this paper and in \cite{kr} we have to use different arguments than for small profinite groups (another reason is for example the fact that invariant subgroups are not necessarily closed and $Acl(\emptyset)$ is not necessarily finite).

We shall need the description of orbits in products of finite, abelian groups from \cite{kr05}. Let $X_{i}$, $i \in \omega$, be a collection of finite abelian groups, and $X =\prod_i X_i$. Let $J_0 \subseteq J_1 \subseteq \dots \subseteq \omega$ be a sequence of finite sets such that $\bigcup J_m =\omega$. Put $X_{J_m}=\prod_{i \in J_m}X_i$. Then, $X$ can be treated as the inverse limit of the $X_{J_m}$. Let $Aut^{0}_J(X)$ be the group of all automorphisms of $X$ respecting this inverse system. For $\eta = (\eta_i: i \in \omega) \in X$, $\eta \! \upharpoonright \! J_m$ will denote the finite sequence $(\eta_i: i \in J_m) \in X_{J_m}$. By \cite[Lemma 4.1]{kr05}, we have the following description of orbits in the profinite group $(X, Aut^0_J(X))$.

\begin{fact}
Let $A$ be a subgroup of $X$ and $\eta=(\eta_1,\dots,\eta_n), \tau=(\tau_1,\dots,\tau_n) \in X^n$. Then, $o(\eta/A)=o(\tau/A)$ iff for all $a \in A$, $k \in {\mathbb Z}$, $l_1,\dots,l_n \in {\mathbb Z}$ and $m \in \omega$ we have
$$k \mid \left( \sum_{i=1}^n l_i\eta_i + a \right)\! \upharpoonright \! J_m \iff k \mid \left( \sum_{i=1}^n l_i\tau_i + a \right)\! \upharpoonright \! J_m.$$
\end{fact}

\begin{corollary}\label{description of orbits}
Assume $X_0=X_1=\dots={\mathbb Z}_p$, where $p$ is a prime number and ${\mathbb Z}_p$ is the cyclic group of order $p$. Let $A$ be a subgroup of $X$, and $\eta=(\eta_1,\dots,\eta_n), \tau=(\tau_1,\dots,\tau_n) \in X^n$. Then, $o(\eta/A)=o(\tau/A)$ iff for all $a \in A$, $l_1,\dots,l_n \in {\mathbb Z}_p$ and $m \in \omega$ we have
$$\left( \sum_{i=1}^n l_i\eta_i + a \right)\! \upharpoonright \! J_m =0 \iff \left( \sum_{i=1}^n l_i\tau_i + a \right)\! \upharpoonright \! J_m = 0.$$
\end{corollary} 

\begin{theorem}\label{counter-example}
For every $\alpha <\omega_1$ there exists a small, $nm$-stable compact $G$-group $(H,G)$ such that $\NM(H)=\alpha$ and $H$ is abelian. Moreover, $H$ has an ascending sequence of $\emptyset$-closed subgroups $H_i$, $i\leq \alpha$, such that $H_i<_{nwd} H_{j}$ whenever $i<j$.
\end{theorem}
{\em Proof.} We can find a descending sequence $(I_i)_{i\leq \alpha}$ of subsets of $\omega$ such that:
\begin{enumerate}
\item $I_i \setminus I_j$ is infinite for all $i<j \leq \alpha$,
\item $I_0 =\omega$ and $I_{\alpha}=\emptyset$.
\end{enumerate}
Put $$H={\mathbb Z}_p^\omega,$$ where $p$ is a prime number.
We equip $H$ with the product topology. Then, $H$ is a profinite group. For $i\leq \alpha$ define $$H_i=\{ \eta \in {\mathbb Z}_p^\omega: \eta(j)=0\; \mbox{for}\; j \in I_i\}.$$
In particular, $H_0=\{0\}$ and $H_\alpha=H$. 
Let $$G=\{ g \in Aut({\mathbb Z}_p^\omega) : g[H_i]=H_i\; \mbox{for every}\; i \leq \alpha \}.$$ 

We will show that $(H,G)$ satisfies the conclusion of the theorem. The following is clear from definitions: $(H,G)$ is a compact $G$-group, $H$ is abelian, all $H_i$'s are $\emptyset$-closed, and $H_i <_{nwd} H_{j}$ whenever $i<j$. So, it remains to show that:    
\begin{enumerate}
\item[(i)] $(H,G)$ is small, 
\item[(ii)] $\NM(H)=\alpha$.
\end{enumerate}

In order to do that, we need to describe orbits on $H$ over finite subsets.
Since $o(a/A)=o(a/\langle A \rangle)$, it is enough to consider orbits over finite subgroups of ${\mathbb Z}_p^\omega$.\\[3mm]
{\bf Claim (Description of orbits)} 
Let $A$ be a finite subgroup of ${\mathbb Z}_p^\omega$, and $\eta, \tau \in{\mathbb Z}_p^\omega$. Define
$$n_\eta=\min \{i\leq \alpha: (\exists a \in A)(\eta + a \in H_i)\}.$$
Then, $o(\eta/A) = o(\tau/A)$ iff $n_\eta=n_\tau$ and $\eta-\tau\in H_{n_\eta}$.\\[3mm]
{\em Proof of Claim.} $(\rightarrow)$ is clear.\\
$(\leftarrow)$. The idea of the proof is to present ${\mathbb Z}_p^\omega$ as the inverse limit of an inverse system of finite groups, find a closed subgroup $G^*$ of $G$ respecting this inverse system, and to show that $\eta$ and $\tau$ are in one orbit under $G^*_A$ by application of Corollary \ref{description of orbits}.

Notice that if $n_\eta=0$, i.e. $\eta \in A$, then $\tau=\eta$, whence $o(\eta/A)=o(\tau/A)$. So, assume $\eta \notin A$. Then, $\tau \notin A$. 

Choose a finite $J_0 \subseteq \omega$ such that for any $a \in A$ we have
$$(\eta + a) \!\upharpoonright \! (J_0 \cap I_{n_\eta -1})  \ne 0\; \mbox{and}\; (\tau + a) \!\upharpoonright \! (J_0 \cap I_{n_\eta -1}) \ne 0.$$
Next, choose any finite $J_1\subseteq J_2 \subseteq \dots \subseteq \omega$ such that $J_0 \subseteq J_1$ and $\bigcup J_m = \omega$. We will show that there exists $f \in Aut^0_J({\mathbb Z}_p^\omega/A,\{H_0\},\{H_1\},\dots)$ (i.e. $f$ is in the pointwise stabilizer of $A$ intersected with setwise stabilizers of $H_0, H_1,\dots$ computed in $Aut^0_J({\mathbb Z}_p^\omega)$) such that $f(\eta)=\tau$. This is enough because $Aut^0_J({\mathbb Z}_p^\omega/A,\{H_0\},\{H_1\},\dots)\le G_A$.

By the choice of $J_m$'s, for all $a \in A$ and $m \in \omega$ we have
$$\left( \eta + a \right)\! \upharpoonright \! J_m \ne 0\; \mbox{and}\; \left( \tau + a \right)\! \upharpoonright \! J_m \ne 0.$$
Thus, by Corollary \ref{description of orbits}, there is $f \in Aut^0_J({\mathbb Z}_p^{\omega}/A)$  such that $f(\eta)=\tau$.

Now, consider any finite $B=\{ b_0,\dots,b_l\} \subseteq H$ with $b_0=0$. We can permute the elements of $B$ so that $b_{j+1} \in H_i$ implies $b_j \in H_i$ for every $i,j$.
For $b \in H$ put $i(b)=\min \{i: b \in H_i \}$. Then, $0=i(b_0) \leq \dots \leq i(b_l)$.

By recursion on $j$, we will construct $f_0,\dots, f_l$ such that:
\begin{enumerate}
\item[1)] $f_j \in Aut^0_J({\mathbb Z}_p^\omega/A)$,
\item[2)] $f_j(\eta)=\tau$,
\item[3)] $f_{j+1} \! \upharpoonright \! \{b_0,\dots,b_j\} = f_j \! \upharpoonright \! \{b_0,\dots,b_j \}$,
\item[4)] $f_{j+1}^{-1} \! \upharpoonright \! \{b_0,\dots,b_j\} = f_j^{-1} \! \upharpoonright \! \{b_0,\dots,b_j \}$,
\item[5)] $f_j(b_j) \in H_{i(b_j)}$,
\item[6)] $f_j^{-1}(b_j) \in H_{i(b_j)}$.
\end{enumerate}

Before the construction, let us explain why this will complete our proof. For any finite $B\subseteq H$ let $G^B$ be the collection of all $f \in Aut^0_J({\mathbb Z}_p^{\omega})$ satisfying 1), 2), 5) and 6) for all elements of $B$. We see that $G^B$ is closed and non-empty. Moreover, for any finite $B_1, \dots, B_n\subseteq H$, $G^{B_1\cup \dots \cup B_n} = G^{B_1} \cap \dots \cap G^{B_n}$, so $G^{B_1} \cap \dots \cap G^{B_n} \ne \emptyset$. By the compactness of $Aut^0_J({\mathbb Z}_p^{\omega})$, the intersection of all such $G_B$'s is non-empty. We also see that any element of this intersection belongs to  $Aut^0_J({\mathbb Z}_p^\omega/A,\{H_0\},\{H_1\},\dots)$ and maps $\eta$ to $\tau$.     

Now, we describe the construction. For the basis step put $f_0=f$. In order to simplify the induction step (i.e.\ to eliminate conditions 4) and 6) from considerations), first we have to show that it is enough to prove the following\\[3mm]
{\bf Subclaim} Let $C \subseteq H$ be finite and $d \in H$ be such that $i(c)\leq i(d)$ for any $c \in C$. Assume $g \in Aut^0_J({\mathbb Z}_p^\omega)$ is such that
\begin{enumerate}
\item[a)] $g \in Aut^0_J({\mathbb Z}_p^\omega/A)$,
\item[b)] $g(\eta)=\tau$,
\item[c)] $i(g(c)) \leq i(d)$ for every $c \in C$.
\end{enumerate}
Then, there is $g_1 \in Aut^0_J({\mathbb Z}_p^\omega)$ satisfying a), b) and 
\begin{enumerate}
\item[d)] $g_1 \! \upharpoonright \! C = g \! \upharpoonright \! C$,
\item[e)] $g_1(d) \in H_{i(d)}$.\\[-2mm]
\end{enumerate}

Let us show that the Subclaim allows us to do the induction step of our  construction. Suppose we have  $f_0,\dots,f_j$ satisfying 1) - 6). First, we apply the Subclaim to $g:=f_j$, $C:=\{b_0,\dots,b_j,f_j^{-1}(b_0),\dots,f_j^{-1}(b_j)\}$ and $d:=b_{j+1}$ in order to get $f_{j+1}':=g_1$ with $f_0,\dots,f_j, f_{j+1}'$ satisfying 1) - 5). Then, we apply the Subclaim (with the reversed roles of $\eta$ and $\tau$) to $g:=f_{j+1}'^{-1}$, $d:= b_{j+1}$ and $C:=\{ b_0,\dots,b_j,f_{j+1}'(b_0),\dots,f_{j+1}'(b_{j+1})\}$ in order to get $f_{j+1}^{-1}:=g_1$ such that $f_0,\dots,f_j,f_{j+1}$ satisfy 1) - 6).\\[3mm] 
{\em Proof of Subclaim.} We can assume $d\ne 0$. There are two cases.\\[2mm]
{\bf Case 1} $i(d)<n_{\eta}$.\\
Since $d\! \upharpoonright \! I_{i(d)} =0$ and for every $c \in C$, $i(d)\geq i(c)$, we have $c \! \upharpoonright \! I_{i(d)} =0$.
Then, by the definition of $n_\eta$, the choice of $J_0$, and the assumption of Case 1,  we see that for every $a \in A$, $k \in {\mathbb Z}_p \setminus \{0\}$ and $c \in \langle C\rangle$ we have
\begin{equation*}\tag{$*$}
\begin{aligned}
(d+k\eta + a+c) \! \upharpoonright \! (J_0 \cap I_{i(d)})&=k(\eta+k^{-1}a)\! \upharpoonright \! (J_0 \cap I_{i(d)})\ne 0,\ \mbox{and}\\
(d+k\tau + a+c) \! \upharpoonright \! (J_0 \cap I_{i(d)})&=k(\tau + k^{-1}a) \! \upharpoonright \! (J_0 \cap I_{i(d)})\ne 0.\end{aligned}\end{equation*}

Consider any $r \in \omega$. We have the following three easy observations (as an example, we give the proof of the third one).
For any $a \in A$ and $c \in \langle C\rangle$ we have: 
\begin{equation*}\tag{$**$}
\mbox{if}\; (d+a+c)\! \upharpoonright \! (J_r \cap I_{i(d)}) \ne 0,\; \mbox{then}\;  (\forall x \in H_{i(d)})((x+a+g(c)) \! \upharpoonright \! J_r \ne 0),
\end{equation*}
\begin{equation*}\tag{$***$}\mbox{if}\; (d+a+c)\! \upharpoonright \! J_{r} \ne 0, \; \mbox{then}\; (g(d)+a+g(c))\! \upharpoonright \! J_{r} \ne 0,
\end{equation*}
\begin{equation*}\tag{$!$}
\mbox{if}\; (d+a+c) \! \upharpoonright \! J_r =0,\; \mbox{then}\; g(d) \! \upharpoonright \! (J_r \cap I_{i(d)}) =0.
\end{equation*}
Indeed, assume $(d+a+c) \! \upharpoonright \! J_r =0$. Then, $(g(d)+a+g(c)) \! \upharpoonright \! J_r =0$. Thus, $(g(d)-d +g(c)-c) \! \upharpoonright \! J_r =0$. But $d\! \upharpoonright \! I_{i(d)} =0$, $c \! \upharpoonright \! I_{i(d)} =0$ and $g(c) \! \upharpoonright \! I_{i(d)} =0$. Therefore, $g(d) \! \upharpoonright \! (J_r \cap I_{i(d)}) =0.$

If for every $r \in \omega$ there are $a \in A$ and $c \in \langle C\rangle$ such that $(d+a+c) \! \upharpoonright \! J_r =0$, then $g(d) \in H_{i(d)}$ by $(!)$, whence $g_1:=g$ works. So, we can assume that there is $r\in \omega$ such that 
\begin{equation*}\tag{$!!$}
(\forall a \in A)(\forall c \in \langle C\rangle)((d+a+c) \! \upharpoonright \! J_r \ne 0),
\end{equation*}
and we choose a minimal $r$ with this property.

Put $J_{-1}=\emptyset$. Define $q:=|J_r \setminus (J_{r-1} \cup I_{i(d)})| \in \omega$.
Using (!!), it is easy to check that ranging over all $a \in A$ and $c \in \langle C\rangle$ such that $(d+a+c) \!\upharpoonright \! (J_r \cap I_{i(d)}) =0$ and $(d+a+c) \! \upharpoonright \! J_{r-1} =0$, we obtain less than $p^q$ values $(a+c) \! \upharpoonright \! J_r$, and therefore less than $p^q$ values $(a+g(c)) \! \upharpoonright \! J_r$. 
Hence, by (!) and the choice of $r$, there exists $d' \in H$ such that:
\begin{itemize}
\item $d' \! \upharpoonright \! I_{i(d)} =0$,
\item $d' \! \upharpoonright \! J_{r-1} = g(d) \! \upharpoonright \! J_{r-1}$,
\item for any $a \in A$ and $c \in \langle C\rangle$ such that $(d+a+c) \!\upharpoonright \! (J_r \cap I_{i(d)}) =0$ and $(d+a+c) \! \upharpoonright \! J_{r-1} =0$ we have $(d'+a+g(c)) \! \upharpoonright \! J_r \ne 0$.
\end{itemize} 

Using this together with $(**)$ and $(***)$, we see that for all $s\geq r$, $a\in A$ and $c \in \langle C\rangle$ one has $(d'+a+g(c)) \! \upharpoonright \! J_s \ne 0$, and hence, by $(!!)$,
$$(d+a+c) \! \upharpoonright \! J_s \ne 0\; \mbox{and}\; (d'+a+g(c)) \! \upharpoonright \! J_s \ne 0.$$
For each $s<r$, $d' \! \upharpoonright \! J_s = g(d) \! \upharpoonright \! J_s$, so we have
$$(d+a+c) \! \upharpoonright \! J_s = 0 \iff (d'+a+g(c)) \! \upharpoonright \! J_s =0.$$

On the other hand, by $(*)$ and the fact that $d' \! \upharpoonright \! I_{i(d)} =0$, we get that for all $a \in A$, $k \in {\mathbb Z}_p \setminus \{0\}$, $c \in \langle C\rangle$ and $s \in \omega$ we have  
\begin{equation*}
(d'+k\eta + a+c) \! \upharpoonright \! J_s \ne 0 \;\, \mbox{and} \;\, (d'+k\tau + a+c) \! \upharpoonright \! J_s \ne 0.
\end{equation*} 

Summarizing, for all $a \in A$, $c \in \langle C\rangle = Lin_{{\mathbb Z}_p}(C)$, $k \in {\mathbb Z}_p$ and $s \in \omega$ we have
$$(d+k\eta + c+a) \! \upharpoonright \! J_s = 0 \iff (d'+k\tau + g(c)+a) \! \upharpoonright \! J_s = 0.$$ 

Using this together with the fact that $g \in Aut^0_J(\mathbb{Z}_p^\omega/A)$ and $g(\eta)=\tau$, we conclude that 
for all $a \in A$, $c \in \langle C\rangle = Lin_{{\mathbb Z}_p}(C)$, $k,l \in {\mathbb Z}_p$ and $s \in \omega$ we have

$$(ld+k\eta + c+a) \! \upharpoonright \! J_s = 0 \iff (ld'+k\tau + g(c)+a) \! \upharpoonright \! J_s = 0.$$

Applying Corollary \ref{description of orbits}, we get the existence of $g_1 \in Aut^0_J({\mathbb Z}_p^\omega/A)$ such that: $g_1(\eta)=\tau$, $g_1 \! \upharpoonright \! C = g \! \upharpoonright \! C$ and $g_1(d)=d' \in H_{i(d)}$. This completes the proof in Case 1.\\[2mm]
{\bf Case 2} $i(d) \geq n_\eta$.\\
The proof in this case is similar, so we only give a sketch. We do not have $(*)$, but the following counterparts of $(**)$, $(***)$ and $(!)$ are present (the assumption of Case 2 is used in the proofs of $(**')$ and $(!')$, which are left to the reader). For any $r \in \omega$,  $a \in A$, $c \in \langle C\rangle$ and $k \in {\mathbb Z}_p$ we have: 
\begin{equation*}\tag{$**'$}
\mbox{if}\; (d+k\eta +a+c)\! \upharpoonright \! (J_r \cap I_{i(d)}) \ne 0,\; \mbox{then}\;  (\forall x \in H_{i(d)})((x+k\tau +a+g(c)) \! \upharpoonright \! J_r \ne 0),
\end{equation*}
\begin{equation*}\tag{$***'$}\mbox{if}\; (d+ k\eta+a+c)\! \upharpoonright \! J_{r} \ne 0, \; \mbox{then}\; (g(d)+k\tau+a+g(c))\! \upharpoonright \! J_{r} \ne 0,
\end{equation*}
\begin{equation*}\tag{$!'$}
\mbox{if}\; (d+k\eta+a+c) \! \upharpoonright \! J_r =0,\; \mbox{then}\; g(d) \! \upharpoonright \! (J_r \cap I_{i(d)}) =0.
\end{equation*}
As in Case 1, using $(!')$, we can assume that there is a minimal $r$ such that
\begin{equation*}\tag{$!!'$}
(\forall a \in A)(\forall c \in \langle C\rangle)(\forall k \in {\mathbb Z}_p)((d+k\eta+a+c) \! \upharpoonright \! J_r \ne 0).
\end{equation*}

Put $J_{-1}=\emptyset$. Define $q:=|J_r \setminus (J_{r-1} \cup I_{i(d)})| \in \omega$.
By $(!!')$, ranging over all $a \in A$, $c \in \langle C\rangle$ and $k \in {\mathbb Z}_p$ such that $(d+k\eta+a+c) \!\upharpoonright \! (J_r \cap I_{i(d)}) =0$ and $(d+k\eta+a+c) \! \upharpoonright \! J_{r-1} =0$, we obtain less than $p^q$ values $(k\tau+a+g(c)) \! \upharpoonright \! J_r$. Hence, by $(!')$ and the choice of $r$, there exists $d' \in H$ such that:
\begin{itemize}
\item $d' \! \upharpoonright \! I_{i(d)} =0$,
\item $d' \! \upharpoonright \! J_{r-1} = g(d) \! \upharpoonright \! J_{r-1}$,
\item for any $a \in A$, $c \in \langle C\rangle$ and $k \in {\mathbb Z}_p$ such that $(d+k\eta+a+c) \!\upharpoonright \! (J_r \cap I_{i(d)}) =0$ and $(d+k\eta+a+c) \! \upharpoonright \! J_{r-1} =0$ we have $(d'+k\tau+a+g(c)) \! \upharpoonright \! J_r \ne 0$.
\end{itemize} 

Using $(**')$, $(***')$, $(!!')$ and the choice of $d'$, we see that 
for all $a \in A$, $c \in \langle C\rangle = Lin_{{\mathbb Z}_p}(C)$, $k \in {\mathbb Z}_p$ and $s \in \omega$ we have
$$(d+k\eta + c+a) \! \upharpoonright \! J_s = 0 \iff (d'+k\tau + g(c)+a) \! \upharpoonright \! J_s = 0.$$

We finish using Corollary \ref{description of orbits} as in Case 1. The proof of the Subclaim and so of the Claim is now completed. \hfill $\square$\\

Now, we will prove (i), i.e.\ smallness of $(H,G)$. Consider any finite $A\subseteq H$. Since we are going to count orbits over $A$, we can assume that $A\le H$. Let $O_1(A):=\{ o(\eta/A): \eta \in H \}$, and $O_1^i(A):=\{ o(\eta/A): \eta \in H, n_\eta=i\}$, $i\leq \alpha$. It is enough to show that each $O_1^i(A)$ is countable.

For $x \in H$ and $i\leq \alpha$ put
$$S_x^i=\{ \eta \in H: \eta \! \upharpoonright \! I_i = x \! \upharpoonright \! I_i\}.$$

By the Claim, the orbits from $O_1^i(A)$ are the sets 
$$S_a^i \setminus \bigcup_{b \in A} S_b^{i-1}, a \in A,$$ 
where for a limit ordinal $\lambda$, $S_b^{\lambda -1}:=\bigcup_{i<\lambda} S_b^i$. 
From this, we see that $O_1^i(A)$ is finite.

It remains to show (ii), i.e. $\NM(H)=\alpha$. Since $H_i$, $i \leq \alpha$, are $\emptyset$-closed and $H_i <_{nwd} H_j$ whenever $i<j$, we easily get $\NM(H) \geq \alpha$.

By induction on $i\leq \alpha$, we will show that for any finite $A<H$ each orbit in $O_1^i(A)$  has $\NM$-rank $\leq i$. This, of course, implies $\NM(H) \leq \alpha$, so the proof of the theorem will be completed.

The basis step is clear as $O_1^0(A)=\{ \{a\}: a \in A\}$. Now, assume that for every finite $A<H$ and $i<j$ each orbit in $O_1^i(A)$ is of $\NM$-rank $\leq i$. Suppose for a contradiction that there is $o(\eta/A) \in O_1^j(A)$ with $\NM(\eta /A)>j$. Then, there is a finite $B>A$ such that $\eta \nmdep_A B$ and $\NM(\eta/B) \geq j$.

On the other hand, we see that each orbit in $O_1^j(A)$ is a dense $G_\delta$ is some $S_a^j$. Thus, since $S_a^j$'s are closed, all 1-orbits over $A$ are non-meager in their relative topologies. In particular, in virtue of \cite[Theorem 2.12]{kr}, $\eta \nmdep_A B$ means that $o(\eta/B) \subseteq_m o(\eta/A)$. Working over $B$, this implies $i:=n_\eta<j$. By induction hypothesis, $\NM(\eta/B)\leq i<j$, a contradiction. \hfill $\blacksquare$

\noindent
Krzysztof Krupi\'nski\\
Instytut Matematyczny, Uniwersytet Wroc\l awski\\
pl. Grunwaldzki 2/4, 50-384 Wroc\l aw, Poland.\\
e-mail: kkrup@math.uni.wroc.pl\\[3mm]
Frank O. Wagner\\
Universit\'e de Lyon; Universit\'e Lyon 1; CNRS\\
Institut Camille Jordan UMR 5208\\ 
B\^atiment Jean Braconnier\\ 
43 boulevard du 11 novembre 1918, F-69622 Villeurbanne-cedex, France.\\
e-mail: wagner@math.univ-lyon1.fr


\begin{thebibliography}{99}

\bibitem{ca}
R. Camerlo, {\em Dendrites as Polish structures}, Proc. Amer. Math. Soc. 139, 2217-2225, 2011.

\bibitem{ew}
D. Evans, F. O. Wagner, {\em Supersimple $\omega$-categorical groups and theories},
J. Symb. Logic 65, 767-776, 2000. 

\bibitem{kr05}
K. Krupi\'nski, {\em Products of finite abelian groups as profinite groups}, J. Alg. 288, 556-582, 2005.

\bibitem{kr05a}
K. Krupi\'nski, {\em Abelian profinite groups}, Fund. Math. 185, 41-59, 2005.

\bibitem{kr}
K. Krupi\'nski, {\em Some model theory of Polish structures}, Trans. Amer. Math. Soc. 362, 3499-3533, 2010.


\bibitem{ne01}
L. Newelski, {\em Small profinite groups}, J. Symb. Logic 66, 859-872, 2001.

\bibitem{ne02}
L. Newelski, {\em Small profinite structures}, Trans. Amer. Math. Soc. 354, 925-943, 2002.

\bibitem{wa03}
F. O. Wagner, {\em Small profinite $m$-stable groups},
Fund. Math. 176, 181-191, 2003.\\[-1mm]

\end{thebibliography}
\end{document}